\documentclass{cmslatex}
\usepackage[paperwidth=7in, paperheight=10in, margin=.875in]{geometry}
 \usepackage[backref,colorlinks,linkcolor=red,anchorcolor=green,citecolor=blue]{hyperref}
\usepackage{amsfonts,amssymb}
\usepackage{amsmath}
\usepackage{graphicx}
\usepackage{cite}
\usepackage{enumerate}
\renewcommand{\theequation}{\arabic{section}.\arabic{equation}}

   \allowdisplaybreaks
\begin{document}
 \title{Second order accurate Dirichlet boundary conditions for linear nonlocal diffusion problems\thanks{Submitted to the editors August 16, 2021. \newline {This research is supported in part  by the NSF DMS-2012562,  DMS-1937254 and  ARO MURI Grant W911NF-15-1-0562.}}}


          \author{Hwi Lee\thanks{Department of Applied Physics and Applied Mathematics, Columbia University, New York, NY
          		10027, USA. {Present Address:  School of Mathematics, Georgia Institute of Technology, Atlanta, GA 30332, USA}. Email:
          		{hl3001@columbia.edu; hlee995@gatech.edu}}
          \and Qiang Du \thanks{Department of Applied Physics and Applied Mathematics, and Data Science Institute, Columbia University, New York, NY
          	10027, USA.   Email:
          	{qd2125@columbia.edu}}}

         \pagestyle{myheadings} \markboth{SECOND ORDER NONLOCAL DIRICHLET BOUNDARY CONDITIONS}{HWI LEE AND QIANG DU} \maketitle

          \begin{abstract}
We present an approach to handle Dirichlet type nonlocal boundary conditions for nonlocal diffusion models with a finite range of nonlocal interactions.  Our approach utilizes a linear extrapolation of prescribed boundary data. A novelty is, instead of using local gradients of the boundary data that are not available a priori, we incorporate nonlocal gradient operators into the formulation to generalize the finite differences-based methods which are pervasive in literature; our particular choice of the nonlocal gradient operators is based on the interplay between a constant kernel function and the geometry of nonlocal interaction neighborhoods. Such an approach can be potentially useful to  address similar issues in peridynamics, smoothed particle hydrodynamics and other nonlocal models. We first show the well-posedness of the newly formulated nonlocal problems and then analyze their asymptotic convergence to the local limit as the nonlocality parameter shrinks to zero. We justify the second order localization rate,  which is the optimal order attainable in the absence of physical boundaries.  
          \end{abstract}

          \section{Introduction}\label{intro}
Nonlocal integro-differential models are increasingly becoming active fields of many scientific and engineering research fronts \cite{bz02,vj12,kc10,dk07,bs16,am17}. Of particular interest to us are those models wherein the nonlocal operators have the feature of a finite range of nonlocal interactions that is parameterized by what is called the nonlocal horizon $\delta$ \cite{dq19,tx13}. On one hand such nonlocal models have been effective in capturing singular physical phenomena such as material defects in the peridynamic setting \cite{ss00}. On the other hand they were the underlying motivations for the development of some meshless numerical methods like the smoothed particle hydrodynamics (SPH) \cite{gr77,ll77}. In this work we are concerned with the following nonlocal linear diffusion models
\[
\mathcal{L}_\delta u(x) := \int_{B_\delta(x)} (u(x)-u(y))w_\delta(x,y)dy = f(x), \quad x \in \Omega
\] on a domain $\Omega$ where $w_\delta(x,y)$ is a nonlocal interaction kernel. Those appear in the bond-based peridynamic models \cite{mt14,ss05} with connections to the SPH setting \cite{dqtx19,dq15}.

There have been a large number of theoretical and numerical studies on the linear nonlocal diffusion models, see for example \cite{dglz13,tn19,mt13,dy17,zwj18,p04,ddggtz20} and the references cited therein. An outstanding issue in the existing body of literature, however,  is to seek suitable nonlocal analogues of local boundary conditions, which in the language of nonlocal vector calculus amounts to imposing the so-called volume constraints \cite{dgr12}. In the presence of physical boundaries an intrinsic challenge is to prescribe the nonlocal interactions in the $\delta$-layer \emph{outside} the domain, which is in fact a distinct feature of generic nonlocal models with a finite $\delta$ range of nonlocal interactions. Our work here provides a nonlocal formulation of volumetric constraints with an emphasis on the convergence rates of the nonlocal solutions to the nonlocal models to those of the corresponding local PDEs as  $\delta \to 0$, when the latter are well-defined. In particular we rigorously establish the second order convergence that are expected to be of optimal order since the nonlocal diffusion operators are second order accurate approximations of the local counterpart in the bulk domain away from the boundaries.

Specifically we focus on the nonlocal Dirichlet volume constraints in parallel with the previous works on the Neumann cases \cite{ttd17,ylty20} wherein the authors rigorously prove the second order convergences of their formulations for one and two dimensional settings, respectively. To our best understanding the predominant approaches for the Dirichlet case have been based on linear extrapolations of the given boundary data on the co-dimension $1$ surface into some volumetric data. The methods of Morris et al \cite{mfz97}, for example, are the ``unofficial standard" \cite{h11} in the SPH context of enforcing no-slip boundary conditions via extrapolation, \textcolor{black}{which can lead to a higher order of accuracy than that of the constant extension approach \cite{magc11} though the latter allows a more straightforward implementation}. However the methods in \cite{mfz97}  have two shortcomings as pointed out in \cite{hp17}; the first is the computational costs of determining distances between particles and the domain boundaries, which could be costly, whereas the second is the implementation of a safety parameter to cap the relative magnitudes of the computed distances for numerical stabilities. The paper \cite{yrs09} proposes a modification of \cite{mfz97}, addressing the complication of non-unique normal distances of a point to the boundary, yet the two drawbacks of the original method remain unresolved in this new approach. The work of \cite{h11} proposes an efficient way to approximate the ratios of distances needed for the extrapolation given in \cite{mfz97}. It, being an estimation after all, might potentially deteriorate the accuracies of the generated numerical solutions. Turning to the peridynamics communities, we note that the recent work \cite{zjcb20} presents an algorithm for nonlocal computations of normal derivatives to locate the mirror nodes in the so-called mirror-based fictions nodes methods, an equivalent to the SPH ghost particle method \cite{cl03}. Although the ghost methods  are designed to avoid the numerical stability issue of  \cite{mfz97}, they could however yield less accurate approximates of the second derivatives.  In \cite{szm21} the authors provide a new version of Taylor series fictitious nodes methods \cite{whzp19} based on the idea of centering the series expansion at the closest node to a given fictitious node. One issue with such an approach, however, is that the constructed approximations of the local derivatives by the conventional finite differences may not be readily accessible depending on the distribution of the nodes. We are hence motivated to consider nonlocal integral approximations of local derivatives which take into account (possibly truncated) circular nonlocal interaction neighborhoods.

Our main contribution is to provide a well-posed nonlocal continuum formulation, the solution to which is proved to be second order convergent to the local counterpart as $\delta \to 0$. What distinguishes our work from the existing literature is that, other than the prescribed local boundary data, we do not require further information about the local solutions. The authors in \cite{dzz19,yyz21} prove the quadratic convergence of their nonlocal solutions given the first order local derivatives of the local solutions on the boundaries.  The latter is generally not given a priori. The strategy adopted in \cite{my21} is to first solve for a family of local problems parameterized by $\delta$ subject to the local boundary data on the outer rim of the $\delta$-layer, and then to use the computed solutions as the volumetric data for the nonlocal problems. Instead of relying on the information about local solutions beyond the given boundary data, we propose to incorporate suitable nonlocal gradient operators into the nonlocal model to mimic the extrapolation of the boundary data to the volumetric data. This amounts in effect to a generalization of the conventional finite-difference based approaches such as in \cite{mfz97} which would correspond to choosing suitable singular delta measures in our nonlocal gradient operators.

We should point out that Zhang and Shi in their recent work \cite{zs21} have already used nonlocal gradient operators to provide $O(\delta^2)$ accurate nonlocal solutions to the nonlocal Poisson problems on manifolds with Dirichlet boundary. They are able to achieve the second order accuracy by enforcing the $O(\delta)$ truncation errors in the inner $\delta$-layer. However we demonstrate that our lower order nonlocal approximations with the $O(1)$ truncation errors in the layer still yield the optimal quadratic convergence rate under suitable assumptions on the nonlocal interaction kernels. This reinforces the view suggested in \cite{dq19} that nonlocal operators can be seen as a continuum weighted average of discrete finite difference operators; the zeroth order truncation errors near the boundaries do not deteriorate the global second order accuracy of the finite difference approximations for the Poisson problems \cite{ym15}. The general form of our nonlocal gradient operators is marked by the extra degree of freedom to specify the nonlocal interaction kernels, which may prove advantageous over the finite difference approximation. For instance we illustrate that one can exploit the interplay between the constant nonlocal kernel and the geometry of nonlocal interaction neighborhoods to implicitly compute the appropriate normalization factors to enforce consistency of the nonlocal operators with the local ones.	The overarching rationale behind our work is to enforce the coupling of the nonlocal diffusion operator in the bulk domain with another \emph{nonlocal} operator for boundary treatment. In order to clearly demonstrate our idea we focus on the one dimensional setting just as the analogous result in one dimensional setting for Neumann boundary conditions \cite{ttd17} is established first before its extension to two dimensions in a subsequent work \cite{ylty20}.

Among other existing works,  the authors of \cite{fry21}  presented explicit rates of localization for the solutions to a vector valued system of nonlocal state-based linear elastic equations subject to general nonlocal Dirichlet-type constraints. Despite the generality of analysis therein, the authors pointed out the sub optimality of their theoretical results which fall $O(\delta^{\frac{1}{2}})$ short of their numerically demonstrated rates. We illustrate here such a theoretical gap can be filled in under a more specialized setting of scalar bond-based models which are indeed later generalized to state-based ones \cite{sew07}. The simplified setting of our work allows us to harness the comparison principles to establish the optimal quadratic convergence rates. At the same time, however, the present work aims to generalize the nonlocal modeling approach behind the conventional approaches \cite{fry21} to impose nonlocal Dirichlet boundary conditions. That is, we allow the material points in the inner $\delta$-layers to take into account possibly different nonlocal Dirichlet boundary values as opposed to constructing \emph{a} set of nonlocal boundary values on the outer $\delta$-layer for \emph{all} the points close to the boundaries. This is achieved by perturbation of the nonlocal diffusion operator by yet another nonlocal operator in the vicinities of the domain boundaries. The present study is in the similar spirit as in various strategies to counteract the surface effects in peridynamics \cite{lb18,msl15} by means of modifying the nonlocal models near the boundaries. We opt for a ``position-aware", in the language of \cite{msl15}, nonlocal formulation in order to align more closely with the SPH standards of Morris et \cite{mfz97} to motivate more practical approaches without requiring the full solution of related local limit.

The paper is organized as follows. In Section \ref{sec:nvp} we present the nonlocal Dirichlet type constrained value problem which incorporates extrapolation via nonlocal gradients. We provide the well-posedness of the equivalent integral equation formulation. In Section \ref{sec:aloc} we prove the weak solution of the nonlocal problem converges to the corresponding local one as the nonlocal horizon vanishes. Moreover we establish the $O(\delta^2)$ convergence rate in uniform norm when the local solution is sufficiently smooth, {followed by the analogous results for some of the existing approaches to enforce the Dirichlet boundary conditions.} In Section \ref{sec:conc} we provide concluding remarks as well as future directions to extend our current work. 

\section{Well-posedness of nonlocal formulation} \label{sec:nvp}
For the rest of the paper we consider the one dimensional domain $\Omega = (0,1)$ unless otherwise noted. 
We seek a nonlocal relaxation of the classical local PDE with homogeneous boundary conditions
\begin{align}
\begin{cases}
-\Delta u(x) = f(x),  & x \in \Omega, \\
u(0) = u(1) = 0,&
\end{cases} \label{eq:le}
\end{align}
in the form of
\begin{align}
\begin{cases}
\widetilde{\mathcal{L}}_\delta u(x):= ({\mathcal{L}}_\delta - {\mathcal{M}}_{\delta}) u(x)  = f(x),  & x \in \Omega, \\
u(x) = 0,  & x \in  (-\delta,0) \cup (1,1+\delta)
\end{cases}.  \label{eq:nlf}
\end{align}
Here $M_\delta$ is a perturbation of $\mathcal{L}_\delta$
\[
\begin{aligned}
\mathcal{M}_{\delta}  u(x) = \int_{B_\delta(x) \backslash \Omega} \left(u(x) + \mathcal{G}^{}_{\delta}u({x})(y-x) \right) w_\delta(x,y) dy
\end{aligned}
\]
where $\mathcal{G}^{}_{\delta}(u)({x})$ is a nonlocal gradient operator given by
\[
\begin{aligned}
&\mathcal{G}^{}_{\delta}u({x}) \\
& = \begin{cases}
\frac{\displaystyle 1}{\displaystyle \int_{B_\delta(x)\cap \Omega} \text{dist}(y,\partial \Omega) \rho (x,y)dy} \begin{cases}
\displaystyle \int_{B_\delta(x)\cap \Omega} u(y) \rho(x,y)dy, & x \in (0,\delta), \\
\displaystyle \int_{B_\delta(x)\cap \Omega} {(-u(y))}\rho(x,y) dy,  & x\in (1-\delta,1),
\end{cases}\\
0, \quad \text{otherwise},
\end{cases}
\end{aligned}
\]
for some nonlocal interaction kernel $\rho$, which we want to choose in order to dispense with calculating the distances $\text{dist}(y,\partial \Omega)$ to the boundaries. Indeed there is such a choice, namely $\rho \equiv 1$, which we utilize for the rest of the paper. Considering without loss of generality $x\in (0,\delta)$, we note that the convexity of $B_\delta(x) \cap \Omega$ gives
\[
\int_{B_\delta(x)\cap \Omega} \text{dist}(y,\partial \Omega) dy = \int_{0}^{x+\delta} |y-0|dy = \int_{0}^{x+\delta} |y-x-\delta|dy = \frac{(x+\delta)^2}{2}.
\] 

Next we specify the assumptions on the kernel $w_\delta(x,y)$ of the nonlocal diffusion operator $\mathcal{L}_\delta$. We assume $w_\delta(x,y) = \frac{1}{\delta^3}w\left(\frac{|x-y|}{\delta}\right)$ where
\begin{equation}
\begin{cases} w \text{ is {continuous},  {nonincreasing} and positive on (0,1) outside which it} \\  \text{ vanishes, and satisfies } \int_\mathbb{R} w(|z|)|z|^2 dz = 2. 
\end{cases} \tag{A1}\label{eq:A1}
\end{equation}
We can then rewrite the equation \eqref{eq:nlf} as the following equivalent integral equation 
\begin{equation}
\widetilde{\mathcal{N}}_\delta  u = f \quad \text{ in } \Omega \label{eq:onl}
\end{equation}
where 
\[
\begin{aligned}
&\widetilde{\mathcal{N}}_\delta u(x) := a_\delta(x) u(x) - \int_\Omega u(y)  \left( w_\delta(x,y) - b_\delta(x) \mathcal{\chi}_{[0,\delta]}(|y-x|)\right)dy. 
\end{aligned}
\]
Here $\mathcal{\chi}_{[\cdot]}$ denotes the characteristic function,
\[
a_\delta(x) =  \int_{\Omega} w_\delta(x,y)dy
\]
and
\[b_\delta(x) = \left \{ \begin{aligned}
& \frac{2}{(x+\delta)^2}  \int_{B_\delta(x)\backslash \Omega} (x-y)w_\delta(x,y)dy,  &  x \in (0,\delta),\\&  \frac{2}{(1-x+\delta)^2}  \int_{B_\delta(x)\backslash \Omega} (y-x) w_\delta(x,y) dy, \quad &  x \in (1-\delta,1), \\
& 0  &  \text{ otherwise}\\
\end{aligned} \right. . 
\]
We note the kernel in the second term of $\widetilde{\mathcal{N}_\delta}$ is in general sign changing and translation-variant. This poses a challenge for us to adapt the technique of \cite{mt13} because the kernel considered therein is translation invariant. Instead we will resort to the idea of Ostrowski's comparison matrix \cite{o37,lpo20} and introduce analogously the comparison operator
\[
\begin{aligned}
&\widetilde{\mathcal{P}}_\delta(u)(x):= a_\delta(x) u(x) - {\int_{0}^{1} u(y) \widetilde{w}_\delta(x,y)  dy}
\end{aligned}
\]
where $\widetilde{w}_\delta(x,y) = |w_\delta(x,y) - b_\delta(x)\chi{1}_{(0,\delta)}(|y-x|)|$. We in turn present the comparison problem 
\begin{equation}
\widetilde{\mathcal{P}}_\delta  v_\delta = f \quad \text{ in } \Omega \label{eq:cnl}.
\end{equation}
The solvability of which will lead to that of the original problem \eqref{eq:onl}. In the meantime, it is not difficult to deduce from the assumptions on $w$ that both $\widetilde{\mathcal{N}}_\delta$ and $\widetilde{\mathcal{P}}_\delta$ are bounded operators on $L^2(\Omega)$, which then leads us to consider $f \in L^2(\Omega)$. With the aforementioned preparations we can now state the well-posedness of \eqref{eq:cnl}.

\begin{proposition} \label{thm:nlcw}
	Suppose that in addition to \eqref{eq:A1}, the kernel $w_\delta$ satisfies
	{
		\begin{equation}
		a_\delta(x) \geq \int_{\Omega}\tilde{w}_\delta(x,y) dy,  \quad  x \in \Omega \tag{A2} \label{eq:A2}
		\end{equation}
		Then $\widetilde{\mathcal{P}}_\delta $ is invertible with  
		\[
		\|(\widetilde{\mathcal{P}}_\delta)^{-1} \|_{2} \leq C
		\]
		for some $C(w,\delta)>0$.}
\end{proposition}

\begin{proof}
	We first argue $\tilde{\mathcal{P}}_\delta$ satisfies the Fredholm alternative. Since $a_\delta(x) > 0$ we have that $a_\delta(x)I$ is an invertible operator. Meanwhile the continuity of $w_\delta$ implies
	$
	\int_\Omega u(y) \tilde{w}_\delta(x,y)dy
	$
	is a compact operator. 
	
	Next we show the kernel of  $\tilde{P}_\delta$ is trivial, which then would imply that $\tilde{P}_\delta$ is invertible. To this end let $v \in L^2(\Omega)$ such that $\tilde{\mathcal{P}}_\delta v \equiv 0$. That is,
	\[
	v(x) = \frac{1}{a_\delta(x)} \int_\Omega v(y) \tilde{w}_\delta(x,y)dy.
	\]  We see $v$ is continuous on $\overline{\Omega}$ due to the continuities of $a_\delta$ and $\tilde{w}_\delta(x,y)$. Then we apply the assumption \eqref{eq:A2} to obtain that $v$ is a constant function. To rule out nonzero constants it is sufficient to show that the inequality in \eqref{eq:A2} is strict in some subdomain of $\Omega$. Indeed for $x \in (0,\delta)$ the non-increasing assumption on $w$ gives 
	\[
	w_\delta(x,y) \geq b_\delta(x)
	\]
	when $ |y-x| \geq x$, so that
	\[
	a_\delta(x) - \int_{0}^{1} \tilde{w}_\delta(x,y) dy
	\geq x b_\delta(x) > 0 \quad x \in \left(\frac{\delta}{2},\frac{3\delta}{4}\right)
	\]
	where the last inequality is due to the positivity of $w_\delta$. 
	
	Finally we conclude the proof by invoking the bounded inverse theorem. 
\end{proof}

\begin{remark}
	Since $w$ is assumed to be a scaled kernel and $\Omega$ is an interval, the condition \eqref{eq:A2} holds for all $\delta \in (0,\delta_0)$ if it is satisfied for some $\delta_0 $. In particular the constant function as well as the piecewise linear function used in \cite{dy19} satisfy \eqref{eq:A2}. An analogues of the condition in the discrete setting would amount to diagonal dominance. 
\end{remark}

Next result we present is concerned with a comparison principle of the comparison problem \eqref{eq:cnl}. One can almost expect the result to hold from the nonnegativity of the kernel $\widetilde{w}_\delta$.  We will utilize the comparison principle in proving the main result of the section, Theorem \ref{thm:nlw} on the well-posedness of the original problem \eqref{eq:onl}.

\begin{lemma}\label{thm:mp}
	Under the same assumptions in Proposition \ref{thm:nlcw}, $\widetilde{\mathcal{P}}_\delta \phi \geq 0$ implies $\phi \geq 0$ for $\phi \in L^2(\Omega)$
\end{lemma} 
\begin{proof}
	A proof is essentially an adaption of the argument in \cite{k79} in the finite dimensional setting, but we include it here for completeness. Let us first write 
	\[
	u = \phi^{+}-\phi^{-}
	\]
	into positive and negative parts. Then since $\tilde{w}_\delta$ is nonnegative we have 
	\[
	\begin{aligned}
	&( \widetilde{\mathcal{P}}_\delta \phi^+,\phi^-)_{L^2} = \underbrace{\int_{0}^{1}a_\delta(x) \phi^{+}(x) \phi^{-}(x) dx}_{=0}  - 		& \int_{0}^{1} \underbrace{\left({\int_{0}^{1} \phi^{+}(y) \widetilde{w}_\delta(x,y)  dy}\right)}_{\geq 0} \phi^{-}(x)dx \leq 0
	\end{aligned}
	\]
	so that if $\phi^{-} \neq 0$
	\[
	0 \leq (\widetilde{\mathcal{P}}_\delta  \phi, \phi^{-}) = (\widetilde{\mathcal{P}}_\delta  \phi^{+}, \phi^{-}) - (\widetilde{\mathcal{P}}_\delta  \phi^{-}, \phi^{-}) <0
	\]
	which is a contradiction. Here the last inequality is due to the fact that $\widetilde{P}_\delta$ has a trivial kernel.
\end{proof}

\begin{theorem} \label{thm:nlw}
	Under the assumptions \eqref{eq:A1} and \eqref{eq:A2}, the same conclusions as in Proposition \ref{thm:nlcw} hold for the problem \eqref{eq:onl} in place of  \eqref{eq:cnl}.
\end{theorem}
\begin{proof}
	By linearity it is sufficient consider the case $f \geq 0$. Since $a_\delta(x)>0$ we can write
	\[
	\widetilde{\mathcal{N}}_\delta = I - \widehat{N}_\delta
	\]
	\[
	\widetilde{\mathcal{P}}_\delta = I - \widehat{P}_\delta
	\]
	where
	\[
	\widehat{N}_\delta \phi(x)= \int_\Omega \phi(y)\frac{w_\delta(x,y) - b_\delta(x) \mathcal{\chi}_{[0,\delta]}(|y-x|)}{a_\delta(x)}dy 
	\]
	\[
	\widehat{P}_\delta \phi(x)= \int_\Omega \phi(y)\frac{|w_\delta(x,y) - b_\delta(x)\mathcal{\chi}_{(0,\delta)}(|y-x|)|}{a_\delta(x)}dy. 
	\]
	By construction we observe
	\[|{\widehat{N}_\delta}^n f(x)| := |\underbrace{\widehat{N}_\delta \circ \dots \circ \widehat{N}_\delta}_{\text{applied n times}} f(x)| \leq {\widehat{P}_\delta}^n f(x), n \in \mathbb{N}, x \in \Omega. \]
	On the other hand we have
	\[
	(\widetilde{P}_\delta)^{-1} f  \geq \sum_{j=0}^{n} {\widehat{P}_\delta}^n f, n \in \mathbb{N}
	\] since we can apply Lemma \ref{thm:mp} to
	\[
	\widetilde{P}_\delta\left((\widetilde{P}_\delta)^{-1} f - \sum_{j=0}^{n} {\widehat{P}_\delta}^j f \right) = {\widehat{P}_\delta}^{n+1} f \geq 0 
	\]
	Consequently we obtain
	\[|  {\sum_{j=0}^{\infty}  {\widehat{N}_\delta}^j f(x)} | \leq (\widetilde{P}_\delta)^{-1} f(x), x \in \Omega.
	\] 
	which proves the existence of $\widetilde{\mathcal{N}}_\delta^{-1}: L^2(\Omega) \to L^2(\Omega)$. As in the proof of Proposition \ref{thm:nlcw}  the proof is complete due to the bounded inverse theorem. 
\end{proof}		

\textcolor{black}{Before we discuss the local limit of our formulation let us briefly mention our nonlocal treatment of the case where the local PDE \eqref{eq:le} is subject to the inhomogeneous Dirichlet boundary conditions  $u(0) = a$ and $u(1) = b$. In that case we may take the linear function $\phi$ which satisfies those boundary conditions, and solve \eqref{eq:onl} with $f$ replaced by $f-\widetilde{\mathcal{N}}_\delta\phi$. Alternatively we may modify the definition of $\mathcal{G}^{}_{\delta}u({x})$ for $x\in (0,\delta) \cup (1-\delta,1)$ into
	\[
	\frac{\displaystyle 1}{\displaystyle \int_{B_\delta(x)\cap \Omega} \text{dist}(y,\partial \Omega) \rho (x,y)dy} \begin{cases}
	\displaystyle \int_{B_\delta(x)\cap \Omega} (u(y)-a) \rho(x,y)dy,  & x \in (0,\delta), \\
	\displaystyle \int_{B_\delta(x)\cap \Omega} {(b-u(y))}\rho(x,y) dy, & x\in (1-\delta,1), 
	\end{cases}
	\] without affecting the validity of our analysis.}

\section{Asymptotic localization as $\delta \to 0$} \label{sec:aloc}   
It is important to analyze consistency between our nonlocal formulation and the corresponding local one as the former is conceived as a relaxation of the latter. We are interested in studying the asymptotic behavior of the nonlocal solution as the nonlocality vanishes, namely its convergence to the local one when the local formulation is mathematically valid. To this end we turn to the variational framework and prove a sharper version of the stability result than Theorem \ref{thm:nlw} wherein the constant $C$ is independent of $\delta$. As in the previous section we will first work with the comparison problem by introducing the quadratic form
\[
\begin{aligned}
Q_\delta(u) := (\widetilde{\mathcal{P}}_\delta  u,u) = \frac{1}{2}  \int_0^{1} \int_0^{1} &(u(y)-u(x))^2 \tilde{w}^{s}_\delta(x,y) dxdy  + \\
& \int_0^{1} u(x)^2 {\left(\int_{0}^{1} (w_\delta(x,y)-\tilde{w}^s_\delta(x,y)) dy\right)} dx. 
\end{aligned}
\] 
One can find in literature the related results when nonlocal interaction kernels are radially symmetric. We would like to make use of those results in our setting and as the first step we show that the symmetric kernel $\widetilde{w}^s_\delta(x,y)$ can be bounded from below by some radially symmetric one. 
\begin{lemma} 
	\label{thm:lb1}
	Suppose $w$ satisfies \eqref{eq:A1}.Then there exists a radial, monotone, non-negative kernel $\rho$ that is compactly supported on $(0,1)$ and strictly positive on $(0,\sigma)$ for some $0 < \sigma$, satisfying
	\[
	\rho_\delta(|x-y|):= \frac{1}{\delta^3}\rho_\delta\left(\frac{|x-y|}{\delta}\right)\leq \widetilde{w}^s_\delta(x,y) \quad \forall x,y \in \Omega 
	\]
	and
	\[
	0 < \int_{\mathbb{R}} \rho_\delta (z)|z|^2 <  \infty.
	\]
\end{lemma}

\begin{proof}
	We consider two cases
	\begin{enumerate}
		\item[(1)] $w$ is constant, i.e. $w_\delta(|x|) = \frac{3}{\delta_3}\mathcal{\chi}_{[0,\delta]}(|x|)$. In this case
		\[
		\widetilde{w}_\delta(x,y) =  \mathcal{\chi}_{[0,\delta]}(|x-y|) \times
		\begin{cases} \frac{6x}{\delta^3(x+\delta)}, &  x\in (0,\delta) \cup (1-\delta,1), \\
		\frac{3}{\delta^3}, & \text{ otherwise}.
		\end{cases}
		\]
		so that
		\[
		\begin{aligned}
		\widetilde{w}^s_\delta(x,y) &= \frac{3}{\delta^3}\mathcal{\chi}_{[0,\delta]}(|x-y|) \times \\
		&\begin{cases}
		\left(\frac{x}{x+\delta} + \frac{y}{y+\delta}\right), & x,y \in (0,\delta) \cup (1-\delta,1), \\
		\left(\frac{x}{x+\delta} + \frac{1}{2}\right), &  x \in (0,\delta) \cup (1-\delta,1), y \in (\delta, 1-\delta),\\
		\left(\frac{y}{y+\delta} + \frac{1}{2}\right), & x \in (\delta, 1-\delta), y \in (0,\delta) \cup (1-\delta,1), \\
		1, & x,y \in (\delta, 1-\delta).
		\end{cases}
		\end{aligned}
		\]
		We set 
		\[
		\rho(|x-y|) = \mathcal{\chi}_{[0,1]}(|x-y|)\frac{3}{2} \frac{|x-y|}{|x-y|+1}.
		\]
		\item[(2)] Otherwise we let
		\[
		\rho(|x|) = \frac{1}{2} \max\left(w(|x|) - 2\int_{0}^{1} y w(|y|) dy ,0\right)
		\]
		which is supported on $(0,\sigma)$ for some $0 < \sigma < 1$ due to \eqref{eq:A1}.
	\end{enumerate}
	In both cases one can verify $\rho_\delta(|x-y|):= \frac{1}{\delta^3}\rho_\delta\left(\frac{|x-y|}{\delta}\right) \leq \frac{1}{2}w_\delta(x,y)$, where $\rho$ satisfies the desired properties, completing the proof. 
\end{proof}

Next result we present is concerned with a sharper version of the inequality in \eqref{eq:A2}. The nonlocal variational framework which we would like to rely on is pertinent to the cased of positive definite nonlocal operators, which we have not yet established for the operator $\widetilde{\mathcal{P}}_\delta$. The following result will help us stride towards that direction.
\begin{lemma}
	\label{thm:lb2}
	Suppose that $w$ satisfies \eqref{eq:A1} and 
	\begin{equation}
	\begin{aligned}
	|\{y \in \Omega: w_\delta(x,y) &- b_\delta(x) \mathcal{\chi}_{[0,\delta]}(|y-x|) \geq 0\}| \\
	& \geq |\{y \in \Omega: w_\delta(x,y) - b_\delta(x) \mathcal{\chi}_{[0,\delta]}(|y-x|) \leq 0\}|,  \quad x \in \Omega.
	\end{aligned}
	\tag{A2s} \label{eq:A2s}
	\end{equation}
	Then 
	\[
	a_\delta(x) - \int_\Omega \tilde{w}_\delta(x,y)dy  \geq  C \delta b_\delta(x), \quad x \in \Omega
	\]
	for some $C>0$ independent of $\delta$.
\end{lemma}

\begin{proof}
	We consider $x \in (0,\delta)$ and recall $w_\delta(|x|) \geq b_\delta(x)$ due to \eqref{eq:A1}. Then we have
	\[
	\begin{aligned}
	\widetilde{\mathcal{P}}_\delta(1) = a_\delta(x) - \int_\Omega \tilde{w}_\delta(x,y)dy &= \int_{0}^{1} \left( w_\delta(x,y) - |w_\delta(x,y) - b_\delta(x)\mathcal{\chi}_{(0,\delta)}(|y-x|)| \right)dy \\
	& = b_\delta(x) (x + (2r(x)-1)\delta) + 2\int_{x + r(x)\delta}^{x+\delta}  w_\delta(x,y)dy 
	\end{aligned}
	\]
	where 
	$$
	r(x) = \sup_{z \in [0,1]} \{ w_\delta(x,x+z\delta) - b_\delta(x) \geq 0 \}.
	$$
	Since $b_\delta(x)$ is decreasing on $(0,\delta)$ it follows from \eqref{eq:A1} and \eqref{eq:A2s} that $r(x)$ is increasing on $(0,\delta)$ with $r(0)\geq\frac{1}{2}$ and $r(\delta) = 1$. Let us consider the different cases respectively.
	\begin{itemize}
		\item If there exists $\lambda \in (0,1)$ such that $\frac{1}{2}< r(\lambda) < 1$, then we have 
		\[
		\begin{aligned}
		\widetilde{\mathcal{P}}_\delta(1) & \geq \begin{cases}
		2\displaystyle \int_{x+r(\lambda)\delta}^{x+\delta}  w_\delta(x,x+y)dy =  \frac{2}{\delta^2} \int_{r(\lambda)}^{1}  w(y) dy, & x \in (0, \lambda \delta), \\
		b_\delta(x)(2r(\lambda) -1)\delta, &  x \in (\lambda \delta, \delta).
		\end{cases}
		\end{aligned}
		\]
		\item Otherwise 
		\[
		r(x) = \begin{cases}
		\frac{1}{2}, \quad& 0 < x < t ,\\
		1, &t < x < 1.
		\end{cases}
		\]
		for some $t \in (0,1)$, so that
		\[
		\begin{aligned}\widetilde{\mathcal{P}}_\delta(1) 
		& \geq \begin{cases}
		2\displaystyle \int_{x+\frac{\delta}{2}}^{x+\delta}  w_\delta(x,x+y)dy = \frac{2}{\delta^2} \int_{\frac{1}{2}}^{1} w(s) ds, & x \in (0, t \delta), \\
		b_\delta(x)\delta \quad &  x \in (t \delta, \delta).
		\end{cases}
		\end{aligned}
		\]
	\end{itemize}
	The assumption \eqref{eq:A1} implies $\int_{s}^{1}  w(y)dy = C_1(s) > 0$ for any $s \in (0,1)$ while we have $\delta b_\delta(0) = \frac{C_2}{\delta^2} \geq \delta b_\delta(x)$, hence proving the desired claim.
\end{proof}

We now present the uniform Poincare's inequality which amounts to variational stabilities of the nonlocal solutions.

\begin{proposition}
	\label{thm:us}
	{Suppose $w_\delta$ satisfies \eqref{eq:A1}, \eqref{eq:A2s} and
		\begin{equation}
		a_\delta(x) \geq \int_{\Omega}\tilde{w}_\delta(y,x) dy  \quad \forall x \in \Omega \tag{A3} \label{eq:A3}.
		\end{equation}
	}  Then there exists $\delta_0 >0$ and $C(w)$ such that
	\[
	\|(\widetilde{\mathcal{N}}_\delta)^{-1}\|_{2} \leq C
	\]
	for $\delta \in (0,\delta_0)$.
\end{proposition}

\begin{proof}
	We observe 
	\[
	|(\widetilde{N}_\delta)^{-1} f(x)| = |(\widetilde{N}_\delta)^{-1} f^{+}(x) -  (\widetilde{N}_\delta)^{-1} f^{-}(x)| \leq 2 (\widetilde{P}_\delta)^{-1} |f|(x), \quad x \in \Omega 
	\]
	from the proof of Theorem \ref{thm:nlw}. 
	Thus it is sufficient to prove 
	\[
	\|(\widetilde{\mathcal{P}}_\delta)^{-1}\|_{2} \leq C
	\]
	for some constant $C>0$ uniformly in $\delta$. We have 
	\[
	\begin{aligned} 
	E_\delta(u) 
	\geq   \frac{1}{2} \int_\Omega\int_\Omega  (u(y)-u(x))^2 & \rho_{\delta}(|x-y|) dxdy +  \\
	&\frac{\tilde{C}}{2}\int_\Omega u^2(x)  \int_{B_\delta(x) \backslash \Omega}\frac{|y-x|}{\delta} w_\delta(|y-x|)dy dx
	\end{aligned} 
	\]
	due to Lemma \ref{thm:lb1},  \eqref{eq:A3} and Lemma \ref{thm:lb2}. If we now let 
	\[
	k_\delta(|x-y|) = \min \left\{ \rho_\delta(|x-y|),\frac{|y-x|}{\delta} w_\delta(|y-x|)\right\}
	\] which has a nonzero finite second moment, then we have
	\[
	\begin{aligned} 
	C E_\delta(u)  \geq  \frac{1}{2} \int_\Omega\int_\Omega (u(y)-u(x))^2 k_{\delta}(|x-y|) dxdy +  \int_\Omega u^2(x) \int_{B_\delta(x) \backslash \Omega} k_\delta(|x-y|) dy dx 
	\end{aligned} 
	\]
	Now we deduce from \eqref{eq:A1} and the monotonity of $\rho_\delta$ that $k_\delta$ satisfies the hypothesis of Theorem 1.3 in \cite{p04} so that we conclude 
	\[
	\begin{aligned} 
	E_\delta(u) \geq  C \|u\|^2_2 \end{aligned}.
	\]	
\end{proof}

\begin{remark}
	The assumption \eqref{eq:A3} would amount to column diagonal dominance in a discrete setting. It is to be checked on $2\delta$-inner layer of $\Omega$, that is, $(0,2\delta) \cup (1-2\delta,1)$, as opposed to the $\delta$ inner layer which is the case for \eqref{eq:A2}. One can verify that it holds true for the constant as well as linear kernels, both of which also satisfy \eqref{eq:A2s}.  
\end{remark}

Another ingredient that is often used to study the local limit is a compactness result such as Lemma 5.2 in \cite{mt14} which is applicable to uniformly (in $\delta$) bounded nonlocal energies. However we are not aware of such result when the kernel is neither radially symmetric nor nonnegative.  Instead we take a detour by considering the following truncation error analysis.

\begin{lemma}
	\label{thm:te}
	For $\phi \in C^4(\bar{\Omega})$ with $\phi(1)=\phi(0)=0$, we have 
	\[
	\begin{aligned}
	&T_\delta(\phi)(x) := (-\Delta + \widetilde{\mathcal{N}}_\delta) \phi(x)  = \begin{cases}
	O(\delta^2), &  x\in (\delta,1-\delta),\\
	O(1) , & \text{ otherwise}.
	\end{cases}.
	\end{aligned}
	\]
	More specifically there exists a constant $M>0$ independent of $\delta$ such that 
	\[
	\begin{aligned}
	&|T_\delta(\phi)(x)| \leq M \delta^3 b_\delta(x)  + O(\delta^2).
	\end{aligned}
	\] 
\end{lemma}

\begin{proof}
	The results are straightforward by a Taylor series expansion and \eqref{eq:A1}.
\end{proof}

Now we establish the convergence in $L^2$ norm of the solutions to the nonlocal problems as the nonlocality parameter vanishes.

\begin{theorem} \label{thm:cl}
	Under the same assumptions as in Proposition \ref{thm:us}, if $u_\delta$ solves \eqref{eq:onl} and $u_0$ solves \eqref{eq:le}, then
	\[
	\lim_{\delta \to 0} \|u_\delta - u_0\|_2  = 0 
	\] 
\end{theorem}
\begin{proof}
	Without loss of generality we consider a sequence of $\delta_n \to 0$ and denote the corresponding nonlocal solution $u_{\delta_n} = (\widetilde{\mathcal{N}}_{\delta_n})^{-1} f$ by $u_n$. Let $\{f_{m}\}_{m=1}^{\infty}$ be a sequence of smooth functions on $\Omega$ such that $\lim_{m\to\infty} \|f_m - f\|_2 = 0$. Now fix $\epsilon > 0$.
	
	Let $u_{m,0}$ solve \eqref{eq:le} with $f_m$ in place of $f$, and let $u_{m,n} = (\widetilde{\mathcal{N}}_{\delta_n})^{-1} f_m$. Then we obtain
	\[
	\begin{aligned}
	&\|u_0-u_n\|_2  \\
	& \leq {\|u_0-u_{{m},0}\|_2} + {\|u_{{m},0}-u_{{m},n}\|_2}+ {\|u_{{m},n}-u_n\|_2} \\
	& \leq  {\|u_0-u_{{m},0}\|_2} + \|(\widetilde{\mathcal{N}}_{\delta_n})^{-1}\|_2 \left(\|T_{\delta_n}u_{{m},0}\|_2 + \|f_{{m}}-f\|_2 \right)  \\
	& \leq {\|u_0-u_{{m},0}\|_2} + C \left(\|T_{\delta_n}u_{{m},0}\|_2 + \|f_{{m}}-f\|_2 \right)
	\end{aligned}
	\]
	where the last inequality is due to Proposition \ref{thm:us}. Now a very crude estimate from Lemma \ref{thm:te} shows
	\[
	\|T_{\delta_n}u_{{m},0}\|_2 \leq  C_1(m) (\delta_n)^{1/2}.
	\]
	Hence if we pick a particular $\hat{m}$ so that
	\[
	{\|u_0-u_{{\hat{m}},0}\|_2} + C  \|f_{\hat{m}}-f\|_2 < \frac{2\epsilon}{3}
	\]
	then we have
	\[
	C \|T_{\delta_n}u_{{m},0}\|_2 \leq C_1(m) (\delta_n)^{1/2} <\frac{\epsilon}{3}
	\]
	for all $n$ sufficiently large, proving the claim.
\end{proof}

We finally turn to the order of convergence rate when the local solution is smooth. We would like to apply the barrier function technique as done in \cite{ttd17}. {At our disposal is the comparison principle Lemma \ref{thm:mp} for the comparison operator $\widetilde{P}_\delta$ so our strategy is to consider a suitable comparison problem. In our explicit construction of a barrier function we exploit the fact that the conclusion of Lemma \ref{thm:lb2} provides a lower bound for $\widetilde{P}_\delta 1$ in terms of $b_\delta$ which is involved in the upper bounds on the  truncation error in Lemma \ref{thm:te}. }

\begin{proposition} \label{thm:1b}
	Suppose $w$ satisfies \eqref{eq:A1} and \eqref{eq:A2s}. Assume further $f$ is regular enough such that the local solution $u$ to \eqref{eq:le} is smooth. If we denote by $u_\delta$ the nonlocal solution to \eqref{eq:onl}, then there exists $\delta_0 >0$ and $C>0$ such that 
	\[
	\|u_\delta - u\|_{\infty} \leq C \delta^2 
	\]  for $\delta \in (0,\delta_0)$
\end{proposition}
\begin{proof}
	We consider 
	\[
	\widetilde{\mathcal{P}}_\delta v_\delta(x) = |T_\delta u_0(x)|
	\]
	so that we have $|e_\delta| \leq 2v_\delta$ as in the proof of Proposition \ref{thm:us}.
	
	Let $\phi(x) = {C} + \psi(x)$ where $ C > 0$ is a constant to be specified, and $\psi$ is the solution to the problem
	\[
	\begin{aligned}
	-\Delta \psi &= 1, \qquad \text{ in } \Omega \\
	\psi &= 0, \qquad \text{ on } \partial \Omega
	\end{aligned}
	\]
	We have then 
	\[
	\begin{aligned}
	&\widetilde{P}_\delta \phi(x) {\geq} {C}(a_\delta(x)-\|\tilde{w}_\delta(x,\cdot)\|_{L^1(\Omega)}) - {C}_1 \left|\int_{\Omega} (x-y) \widetilde{w}_\delta(x,y)dy \right| \\
	& {-\psi^{\prime\prime}(x)}\int_{\Omega} \frac{(y-x)^2}{2} \widetilde{w}_\delta(x,y) dy    \\
	& \geq  C C_2 \delta b_\delta(x) - {C}_1 \left|\int_{\Omega} (x-y) \widetilde{w}_\delta(x,y)dy \right|  + \int_{\Omega} \frac{(y-x)^2}{2} \widetilde{w}_\delta(x,y) dy\\
	&\geq C {C_2\delta} b_\delta(x) - C_1 \left(\left|\int_{\Omega} (x-y) {w}_\delta(x,y)dy\right|  + b_\delta(x) \underbrace{\left|\int_\Omega (x-y)  \mathcal{\chi}_{(0,\delta)}(|y-x|) dy \right|}_{{O(\delta^2)}} \right) \\ 
	&\int_{\Omega} \frac{|y-x|^2}{2} {w}_\delta(x,y) dy    - b_\delta(x) \underbrace{\left( \int_{\Omega} \frac{|y-x|^2}{2} \mathcal{\chi}_{(0,\delta)}(|y-x|) dy \right)}_{{O(\delta^3)}} 
	\end{aligned}
	\]
	for some constants $C_1, C_2 > 0$. Here the first, second and third inequalities hold due to the non-negativity of $\psi$, Lemma \ref{thm:lb2}, and the property $\widetilde{w}(x,y)= \widetilde{w}(x,z)$ for $|y-x| = |z-x|$, respectively. Thus for sufficiently large $\tilde{C}>0$ we can write  
	\[
	\begin{aligned}
	\widetilde{P}_\delta \phi(x) \geq M\delta b_\delta(x) - C_1 \left|\int_{\Omega} (x-y) {w}_\delta(x,y)dy\right|  + \int_{\Omega} \frac{|y-x|^2}{2} {w}_\delta(x,y) dy
	\end{aligned}
	\]
	{But then since $\left|\int_{\Omega} (x-y) {w}_\delta(x,y)dy\right| \geq \delta^2 b_\delta(x)$} there exists $\delta_0>0$ and $C>0$ such that
	\[
	\begin{aligned}
	\widetilde{P}_\delta \phi(x) \geq {M \delta b_\delta(x) } + \int_{\Omega} \frac{|y-x|^2}{2} {w}_\delta(x,y) dy \geq {M \delta b_\delta(x) + C }
	\end{aligned}
	\]
	for all $\delta < \delta_0$,  {where the assumption \eqref{eq:A1} is used in the last inequality. Therefore we obtain
		\[
		\begin{aligned}
		|v_\delta(x)| \leq \left(\sup_{x\in \Omega} \phi(x)\right) \left(\sup_{x\in \Omega}\frac{|T_\delta u_0(x)|}{\widetilde{\mathcal{P}}_\delta \phi(x)} \right) \leq C_1 \sup_{x\in \Omega}\frac{M_1 \delta^3b_\delta(x)+O(\delta^2)}{M \delta b_\delta(x)+C} \leq  C \delta^2
		\end{aligned}
		\]
		where the first and second inequalities are due to Lemma \ref{thm:mp} and Lemma \ref{thm:te}, respectively, completing the proof. 
	}
\end{proof}
{
	\begin{remark}
		In comparison with Theorem \ref{thm:cl} we have removed the assumption \eqref{eq:A2s} at the expense of a stronger regularity assumption on $f$.
\end{remark}}

Based on our choice of the barrier function in the proof of Proposition \ref{thm:1b} we can in fact revisit the existing methods of enforcing nonlocal Dirichlet boundary conditions to prove their convergence rates as $\delta \to 0$. The two approaches that we focus on are the constant extension methods \cite{magc11} and the methods of Morris et al \cite{mfz97}. Notwithstanding their numerically demonstrated performances, we are not aware, to our best understanding, of existing theoretical justification for the convergence rates of their associated nonlocal \emph{solutions}. 

We first turn to the constant extension methods wherein the nonlocal volumetric constraints are prescribed via constant extension of the local boundary values, that is,  $u(x) = 0$ for $x \in (-\delta,0) \cup  (1,1+\delta)$. Hence the nonlocal continuum formulation, in the case of homogeneous boundary conditions as before, is given by
\begin{align}
\begin{cases}
{\mathcal{L}}_\delta u(x)  = f(x),  \quad & x \in \Omega, \\
u(x) = 0, \quad  &x \in (-\delta,0) \cup (1,1+\delta).
\end{cases}. \label{eq:cem}
\end{align}
One can find a multitude of thorough mathematical investigations on this formulation; see, for example, \cite{dq19} and the references cited therein. In particular, the recent work \cite{yyz21} has established the first order localization rate of the formulation. Nevertheless an alternative proof of the linear convergence can be established if we note that the truncation errors of the formulation, in the sense of Lemma \ref{thm:te}, are bounded by ${M \delta^2 b_\delta(x)}$ for some $M$ independent of $\delta$ on $(-\delta,0) \cup (1,1+\delta)$. With the same choice of the barrier function as before we can prove
\begin{proposition}
	Suppose $w$ satisfies \eqref{eq:A1}. Assume further $f$ is regular enough such that the local solution $u$ to \eqref{eq:le} is smooth. If we denote by $u_\delta$ the nonlocal solution to \eqref{eq:cem}, then there exists $\delta_0 >0$ and $C>0$ such that 
	\[
	\|u_\delta - u\|_{\infty} \leq C \delta
	\]  for $\delta \in (0,\delta_0)$.
\end{proposition}

Next we analyze the methods of Morris et al \cite{mfz97}, which is based on extrapolation via finite difference approximations of local gradients. The underlying continuum nonlocal formulation amounts to
\begin{align}
\begin{cases}
{\mathcal{L}}_{\delta,M} u(x):= ({\mathcal{L}}_\delta + {\mathcal{F}}_{\delta}) u(x)  = f(x),  \quad & x \in \Omega ,\\
u(x) = 0, \quad & x \in  (-\delta,0) \cup (1,1+\delta).
\end{cases}
\end{align}
where
\[
\begin{aligned}
&\mathcal{F}^{}_{\delta}(u)({x})  \\
& = \begin{cases}
u(x)
\displaystyle \int_{B_\delta(x)\backslash \Omega} \left| (p(x)-y) \cdot \frac{x-p(x)}{|x-p(x)|^2}\right|
w_\delta(x,y)dy, \ & x \in (0,\delta) \cup (1-\delta,1),\\
0, \quad & \text{otherwise}.
\end{cases}.
\end{aligned}
\]
Here $p(x)$ denotes the closest point to $x$ on $\partial \Omega$.  As in our formulation we can rewrite the formulation into

\begin{equation}
{\mathcal{N}}_{\delta,M} u = f \quad \text{ in } \Omega \label{eq:mtf}
\end{equation}
where ${\mathcal{N}}_{\delta,M}$ is given by
\begin{equation*}
\begin{aligned}
{\mathcal{N}}_{\delta,M} u(x) &:= \left({\int_{B_\delta(x)\backslash \Omega} w_\delta(|x-y|)dy}\right)u(x) + \mathcal{F}^{}_{\delta}(u)({x})  \\
& + \int_\Omega (u(x)-u(y))   w_\delta(|x-y|) dy.
\end{aligned} 
\end{equation*} 
We then define the associated energy space $V_\delta = \{u\in L^2(\Omega): ({\mathcal{N}}_{\delta,M} u,u) < \infty \}$  for which we have the following characterization. 
\begin{lemma}  \label{thm:mnfs}
	Assume $w$ satisfies \eqref{eq:A1}. Then \textcolor{black}{the problem \eqref{eq:mtf} is well-posed over the space $V_\delta$, a weighted $L^2$ space with the weight $\eta_\delta$ 
		\[
		\eta_\delta(x) = 1 + {\int_{B_\delta(x)\backslash \Omega}}\frac{|(p(x)-y)\cdot (x-p(x))| }{\left|x-p(x)\right|^2}w_\delta(x,y)dy.
		\]} Moreover $H^1_{0}(\Omega)$ is a proper subset of $V_\delta$.
\end{lemma}
\begin{proof}
	\textcolor{black}{Since the operator 
		$
		\int_\Omega (u(x)-u(y))   w_\delta(|x-y|) dy
		$
		is a bounded operator on $L^2(\Omega)$ we see that $V_\delta$ is the weighted $L^2$ space which is complete, hence the well-posedness follows due to the Riesz representation theorem. For the second claim we consider}
	
	\[
	\begin{aligned}
	& \int_{0}^{\delta} \frac{u^2(x)}{x} \left( \int_{x-\delta}^{0} |y| w_\delta(x,y)dy \right) dx   \\
	& \leq \left(\int_\mathbb{R} |z|w_\delta(|z|) dz \right) \left( \int_0^{\delta} \frac{u^2(x)}{x^2}dx  \right) \leq \left(\int_\mathbb{R} |z|w_\delta(|z|) dz \right) \left(\int_{0}^{\delta} (u^\prime(x))^2  \right)
	\end{aligned}
	\] 
	for $u \in H^1_0(0,1)$, where the last inequality is due to the Hardy's inequality \cite{bm97}. Finally we can take $\phi(x) = \sqrt{x}-x$ to show the inclusion is proper, completing the proof. 
\end{proof}

{The inclusion of the Sobolev space $H^1_0(\Omega)$ in the nonlocal space $V_\delta$ clarifies the connection that the integral equation \eqref{eq:mtf} is a suitable nonlocal candidate to approximate the local differential equation \eqref{eq:le}. At the same time we point out $V_\delta$ is strictly included in the solution space $L^2(\Omega)$ of our formulation \eqref{eq:onl} since a nonzero constant function does not belong to $V_\delta$; the two may indeed be viewed as exploiting distinct degrees of nonlocal relaxations.} A simple calculation shows that the same truncation error estimates hold for ${\mathcal{N}}_{\delta,M}$ as in Lemma \ref{thm:te}, and the same barrier function as before can be used to show the following.

\begin{proposition}
	Suppose $w$ satisfies \eqref{eq:A1}. Assume further $f$ is regular enough such that the local solution $u$ to \eqref{eq:le} is smooth. If we denote by $u_\delta$ the nonlocal solution to \eqref{eq:mtf}, then there exists $\delta_0 >0$ and $C>0$ such that 
	\[
	\|u_\delta - u\|_{\infty} \leq C \delta^2
	\]  for $\delta \in (0,\delta_0)$.
\end{proposition}

In comparison with our approach, the trivial extension  \cite{magc11} and the linear extrapolation strategies in \cite{mfz97} admit a wider range of nonlocal interaction kernels due to their more relaxed assumptions on the kernels. The simplicity of the constant extension method is its attractive feature which comes with the price of being only first order accurate. As opposed to the method of Morris et al \cite{mfz97}, our formulation is based on the nonlocal operator that is not self-adjoint. Our approach however is closely related to the work \cite{zs21} by Zhang and Shi which utilizes yet another non self-adjoint nonlocal operator in order to obtain the second order accurate approximation. In contrast with our choice of $w_\delta$ they adopt more regularly scaled kernel $W_\delta(x,y) = C_\delta W\left( \frac{|x-y|^2} {4\delta^2}\right)$, where $W$ is radial, non-negative, smooth and compactly supported on $(-2\delta,2\delta)$, and $C_\delta$ is a normalization constant. Their formulation can be expressed in the current $1D$ setting as 
\begin{equation}
\begin{aligned}
&\int_\Omega (u(x)-u(y)) \frac{W_\delta(|x-y|)}{\delta^2}dy -2 \widehat{\mathcal{G}}_\delta u(p(x)) \overline{W}_\delta(|x-p(x)|) \\
& = \int_\Omega f(y) \overline{W}_\delta(|x-y|)dy +  f(p(x))|x-p(x)| \overline{W}_\delta(|x-p(x)|)
\end{aligned} \label{eq:szf}
\end{equation}
where $\overline{W}_\delta(|x-y|) = C_\delta \overline{W}(\frac{|x-y|}{4\delta^2})$, $\overline{W}(r) = \int_r^{1} W(z)dz$, and
\[
\begin{aligned}
&\widehat{\mathcal{G}}_{\delta}u({x}) \\
& = \begin{cases}
\frac{-1}{2\delta^2 \overline{\overline{W}}_\delta(0)} 
{\displaystyle
	\int_{0}^{2\delta}  \left(  u(y)  \overline{W}_\delta(|y|) + \delta^2  f(y)  \overline{\overline{W}}_\delta(|y|)\right) dy}, & x \in (0,2\delta),\\
\frac{-1}{2\delta^2 \overline{\overline{W}}_\delta(0)} 
{\displaystyle \int_{1-2\delta}^{1}  \left(  u(y)  \overline{W}_\delta(|1-y|) + \delta^2  f(y)  \overline{\overline{W}}_\delta(|1-y|)\right) dy}, & x \in (1-2\delta,1),\\
0, \quad  &\text{otherwise}.
\end{cases}
\end{aligned}
\]
for $\overline{\overline{W}}_\delta(|x-y|) = C_\delta \overline{\overline{W}}(\frac{|x-y|}{4\delta^2}), \overline{\overline{W}}(r) = \int_r^{1} \overline{W}(z)dz$. The operators $\widehat{\mathcal{G}}_{\delta}$ approximate the outward normal local derivatives with $O(\delta^2)$ error, contributing to $O(\delta)$ truncation error of the overall formulation in $(0,2\delta) \cup (1-2\delta,1)$, which is one $\delta$ order higher than our formulation. Nevertheless we are still able to maintain the second order accuracy without the dependence of our nonlocal gradient operator  ${\mathcal{G}}_{\delta}$ on the source term $f$ or the modification of the right hand side in the nonlocal equation.\eqref{eq:onl}. \textcolor{black}{Our approach also admits a wider class of nonlocal kernels for the nonlocal gradient operators ${\mathcal{G}}_{\delta}$ including those that are constructed from the kernels of the nonlocal diffusion operators as done in Zhang and Shi's formulations. This may prove advantageous in terms of enhanced modeling capabilities. Apart from the different scalings of the kernels their approach can be seen as a special case of ours in the sense that the nonlocally extrapolated boundary values in the outer $\delta$ collars are no longer position dependent in the interior $\delta$-layers.}

Notwithstanding the lack of self adjointness in our approach we illustrate in the sequel a scenario in two dimensions where the method of Morris et al \cite{mfz97} may lose its optimal rate of convergence in the presence of a relatively simple circular boundary, for which we propose ours as a viable alternative. We remark that our approach is based on volume integrals only, as opposed to the 2D formulation of Zhang and Shi  \cite{zs21} wherein the boundary surface (curve) integrals are also computed.

\subsection{A two dimensional case study}
In this subsection we present an extension of our formulation to a two dimensional setting focusing on its quadratic convergence rate. We promote our extension as a more accurate approach than the method of Morris et al \cite{mfz97}, which is expected to provide the sub-optimal convergence rate in our chosen two dimensional domain. Specifically let us consider as in \cite{cdprx20} the punctured periodic domain $\Omega = [-2,2]^2 \backslash U $ where $U  = \{ \vec{x} = (x_1,x_2) \in \mathbb{R}^2: |\vec{x} \leq 1\}$. We assume a sufficiently smooth $f$ such that there exists a smooth  $u$ which is periodic on $\partial \Omega \backslash \partial U$ and
satisfies
\begin{align}
\begin{cases}
-\Delta u(x) = f,  &\quad x \in \Omega, \\
u(x) = 0, &\quad x \in \partial U.\\
\end{cases}  \label{eq:2dle}
\end{align}
We choose for simplicity the constant kernel $w_\delta(\vec{x},\vec{y}) = \frac{16}{\delta^4}\chi_{(0,\delta)}(|\vec{x}-\vec{y}|)$ which satisfies the normalization condition $\int_{\mathbb{R}^2} w_\delta(|\vec{x}|)|\vec{x}|^2 d\vec{x} = 4$. For concreteness we fix $\vec{x}^{\star} = (1+\epsilon,0)$ where  $0<\epsilon < \delta$. We first recall the nonlocal operator in our $1D$ formulation which can be rewritten as
\[
\begin{aligned}
\int_0^1 (u(x)-u(y)) &w_\delta(x,y)dy + \left({\int_{0}^{1} \min(|z|,|z-1|) \chi_{(0,\delta)}(|x-z|)dz}\right)^{-1} \cdot \\ 
&  \left({\int_{0}^{1} u(z) \chi_{(0,\delta)}(|x-z|)dz}\right) \left|\int_{0}^{1} (y-x) w_\delta(x,y)dy \right|.
\end{aligned}
\]
We benchmark this in our current two dimensional setting by first defining the direction vector
\[
n_\delta(\vec{x}^{\star}) = \frac{\frac{16}{\delta^4}{\displaystyle \int_{\Omega} (\vec{y}-\vec{x}^\star) \chi_{(0,\delta)}(|\vec{y}-\vec{x}^\star|)d\vec{y}}}{\left| \frac{16}{\delta^4}{\displaystyle  \int_{\Omega} (\vec{y}-\vec{x}^\star) \chi_{(0,\delta)}(|\vec{y}-\vec{x}^\star|) d\vec{y} } \right|} = (1,0).
\]
Next we compute a nonlocal directional derivative in the direction of $n_\delta(\vec{x}^\star)$ 
\[
\mathcal{G}_\delta^{n_\delta(\vec{x}^\star)}u(\vec{x}^\star) = \left(\int_{\Omega} |q_{n_\delta(\vec{x}^\star)}(\vec{y})-\vec{y}| \chi_{\mathcal{I}_\delta(\vec{x}^\star)}(\vec{y}) d\vec{y}
\right)^{-1} \int_{\Omega} u(\vec{y}) \chi_{\mathcal{I}_\delta(\vec{x}^\star)}(\vec{y}) d\vec{y}
\]
where $q_{n_\delta(\vec{x}^\star)}(\vec{y})$ is the projection of $\vec{y}$ onto $\partial \Omega$ along $n_\delta(\vec{x}^\star)$ and $\mathcal{I}_\delta(\vec{x}^\star) = B_\delta(\vec{x}^\star) \cup \{\vec{y} \in \mathbb{R}^2: \vec{y} \cdot n_\delta(\vec{x}^{\star}) = y_1 \leq 0 \}$.  Now the geometry of $\mathcal{I}_\delta(\vec{x}^\star) \cap \Omega$ allows us to bypass the calculations of $q_{n_\delta(\vec{x}^\star)}(\vec{y})$ by using the distance along $n_\delta(\vec{x}^\star)$ to the segment  $\{ \vec{y} \in \partial B_\delta(\vec{x}^\star) \cap \Omega: (\vec{y}-\vec{x}^\star)\cdot \vec{n}_\delta(\vec{x}^\star) \geq 0 \}$ 
\[
\int_{\Omega} |q_{n_\delta(\vec{x}^\star)}(\vec{y})-\vec{y}| \chi_{\mathcal{I}_\delta(\vec{x}^\star)}(\vec{y}) d\vec{y} = \int_\Omega (1+\epsilon+\sqrt{\delta^2-y_2^2}-y_1)\chi_{\mathcal{I}_\delta(\vec{x}^\star)}(\vec{y})  dy_1dy_2.
\]
Finally our nonlocal formulation at $\vec{x} = \vec{x}^\star$ is given by
\[
\int_{\Omega} (u(\vec{x}^{\star})-u(\vec{y})) \frac{16}{\delta^4}d \vec{y} + \mathcal{G}_\delta^{n_\delta(\vec{x}^\star)}u(\vec{x}^\star) \int_{\Omega} (\vec{y}-\vec{x})\cdot n_\delta(\vec{x}^\star) \frac{16}{\delta^4} d\vec{y} = f(\vec{x}^\star)
\]
which has the $O(1)$ truncation error by construction. As in our analysis of the $1D$ case we can rewrite our formulation into the integral equation
\begin{equation}
a_\delta(\vec{x}) u(\vec{x}) - \int_\Omega u(\vec{y})  \left( \frac{16}{\delta^4} - \widehat{b}_\delta(\vec{x}) \mathcal{\chi}_{\mathcal{I}_\delta(\vec{x})}(\vec{y})\right)d\vec{y} = f(\vec{x}) \label{eq:2dne}
\end{equation}
where 
\[
\widehat{b}_\delta(\vec{x}) = 
\left(\int_{\Omega} |q_{n_\delta(\vec{x})}(\vec{y})-\vec{y}| \chi_{\mathcal{I}_\delta(\vec{x})}(\vec{y}) d\vec{y}
\right)^{-1} \left(\int_{\Omega} (\vec{y}-\vec{x})\cdot n_\delta(\vec{x}) \frac{16}{\delta^4} d\vec{y} \right)
\]
if $\vec{x}$ belongs to the $\delta$-layer $\{\vec{z}\in \Omega: \text{dist}(\vec{z},\partial U) < \delta\}$ whereas $\mathcal{L}_\delta u(\vec{x}) = f$ otherwise. Then we have the following quadratic convergence result.

\begin{proposition}
	The nonlocal equation \eqref{eq:2dne} is well-posed over $L^2(\Omega)$. Moreover its unique solution $u_\delta$ converges in $L^\infty$ to the local solution $u_0$ of \eqref{eq:2dle} at the rate of $O(\delta^2)$ as $\delta \to 0$.
\end{proposition}
\begin{proof}
	The similar line of analysis as in $1D$ can be adopted, hence we only sketch the key steps without detailed calculations. The first is to check that the comparison operator satisfies
	\[
	a_\delta(\vec{x}) - \int_\Omega \left| \frac{16}{\delta^4} - \widehat{b}_\delta(\vec{x}) \mathcal{\chi}_{\mathcal{I}_\delta(\vec{x})}(\vec{y})\right|d\vec{y} \geq C\delta \int_{\Omega} (\vec{y}-\vec{x})\cdot n_\delta(\vec{x}) \frac{16}{\delta^4} d\vec{y}
	\]
	for some $C > 0$ in parallel with Lemma \ref{thm:lb2}. The second is to apply the Fredholm argument as in Theorem \ref{thm:nlw} for the well-posedness of the comparison problem, hence the original nonlocal equation \eqref{eq:2dne}. Lastly we use the comparison principle together with the barrier function $\tilde{C}+\phi$ for sufficiently large $\tilde{C} > 0$, where $\phi$ solves \eqref{eq:2dle} with $f \equiv 1$.
\end{proof}

Returning to the method of Morris et al  \cite{mfz97}, let us note $p(\vec{x}^{\star}) = (1,0)$ and the leading term in the truncation error $T_\delta(\vec{x}^{\star})$ is  given by 
\[ 
\frac{16}{\delta^4}\partial_{1} u(1+\epsilon,0)\underbrace{\left(\int_{B_\delta(\vec{x}^{\star})\cap \Omega} (y_1-1-\epsilon) dy_1dy_2 + \int_{B_\delta(\vec{x}^{\star})- \Omega} (y_1-1) dy_1dy_2\right)}_{O(\delta^3)} = O\left(\frac{1}{\delta}\right).
\] 
It should be noted that this estimate is sharp since the term in the parenthesis is  for instance, $\approx 2.9749e-04,  3.7654e-05,   4.7393e-06,   5.9457e-07$ when $\epsilon = \frac{\epsilon}{3}$ and $\delta = 0.1,0,05,0.025,0.00125$, respectively. Consequently there is a loss of opitmal second order convergence which can be seen from the equation for the error $\widehat{\mathcal{N}}_{\delta} e_\delta = T_\delta$ at $\vec{x} = \vec{x}^\star$, that is 
\[
\begin{aligned}
\frac{\delta^4}{16} & \underbrace{\left({\int_{B_\delta(\vec{x}^{\star}) \cap \Omega}} d\vec{y} + \frac{1}{\epsilon^2}\int_{B_\delta(\vec{x}^{\star})\backslash \Omega} \left| (p(\vec{x}^{\star})-\vec{y}) \cdot (\vec{x}^{\star}-p(\vec{x}^{\star}))\right| d\vec{y} \right)}_{O\left(\frac{1}{\delta^2}\right)}  e_\delta(\vec{x}^{\star}) \\
&+ \frac{C}{\delta^4} \int_{\Omega \cap B_\delta(\vec{x}^{\star})}  e_\delta(\vec{y}) d\vec{y}  =  T_\delta(\vec{x}^{\star})
\end{aligned}
\]
This shows that if $e_\delta = O(\delta^2)$ then the  hand side would be $O(\delta)$, which contradicts the right hand side being $O(\frac{1}{\delta})$. Our analysis in Section \ref{sec:aloc} suggests the action of $\widetilde{\mathcal{N}}_{\delta}^{-1}$ lifts the order of the truncation error by $O(\delta^2)$ in the $\delta$-layer, hence we can instead expect to get only $O(\delta)$ error in the uniform norm. 

\section{Conclusion}
\label{sec:conc}
In this work we have presented a nonlocal formulation to enforce local Dirichlet type boundary conditions in the context of nonlocal linear diffusion problems. The key ingredient of our method is the application of nonlocal gradient operators for a linear extrapolation under  suitable conditions on nonlocal interactions kernels to ensure the well-posedness of the resulting nonlocal formulations. Of perhaps more interests to scientific communities at large is our justification of the second order rate in $\delta$ at which the solution to the nonlocal problem converges the local counterpart uniformly in $\Omega$. We point out that our convergence analysis supplements the previous work \cite{magc11} which focuses on the consistency between the nonlocal and local operators as opposed to the solutions of their continuum formulations.

Part of our ongoing research efforts is to extend our current work to more general two dimensional domains for which characterizations of their geometries need be carefully taken into account. Future investigations will also include theoretical analysis of vector valued systems of nonlocal equations, possibly complementing the work \cite{fry21}, so that they can be applied to address similar issues in the study of peridynamics, smoothed particle hydrodynamics (SPH) and other nonlocal models. In particular we are interested in the viscosity formulations of the velocities subject to the no-slip boundary conditions in the SPH-like methods given their popularity as practical simulation tools. This in turn is likely to warrant thorough numerical studies of our nonlocal boundary formulations with respect to the existing criteria for effective numerical evaluations of nonlocal models, such as the notion of asymptotic compatibilities \cite{td14}.

\section*{Acknowledgments}
The authors would like to thank the members of the Computational
Mathematics and Multiscale Modeling (CM3) group at Columbia University for discussions.

\bibliographystyle{siamplain}
\bibliography{cms_manuscript_references}

          \end{document}